\documentclass[12pt]{article}

\usepackage{amsmath,amssymb,latexsym}
\usepackage{amsfonts,amstext,amsthm,amscd}
\usepackage{graphicx}
\usepackage{enumerate}
\usepackage{amsthm}
\theoremstyle{plain}

\usepackage{color}

\newtheorem{theorem}{Theorem}[section]

\newtheorem{teo}[theorem]{Theorem}

\def\R{\mathbb{R}}

\def\R{\mathbb{R}}

\def\E{\mathbb{E}}

\def\N{\mathbb{N}}
\def\P{\mathbb{P}}
\def\F{\mathcal{F}}

\newcommand{\ter}{\hspace{\stretch{3}}$\square$\\[1.8ex]}

\def\l{\langle}

\def\r{\rangle}

\title{Super-Brownian motion: $L^p$-convergence of martingales through the pathwise spine decomposition}
\author{ A.E. Kyprianou\footnote{{\sc University of Bath, UK.} E-mail: a.kyprianou@bath.ac.uk}\ \ \ \  and \ \ \ \ 
 A. Murillo-Salas\footnote{{\sc  Universidad de Guanajuato, M\'exico.} E-mail: amurillos@ugto.mx} \\
  {\it On the occasion of the 60th birthday of \textcolor{black}{Sergei} Kuznetsov}
 }
\begin{document}


\maketitle
\begin{abstract}\noindent 
Evans \cite{Evans} described the semi-group of a superprocess with quadratic branching mechanism under a martingale change of measure in terms of the semi-group of an immortal particle 
and the semigroup of the superprocess prior to the change of measure. This result, commonly referred to as the spine decomposition, alludes to a 
pathwise decomposition in which independent copies of the original process `immigrate' along the path of the immortal particle. For branching particle 
diffusions the analogue of this decomposition has already been demonstrated in the pathwise sense, see for example \cite{Kyp, HH2009}. The purpose of
 this short note is to exemplify a new {\it pathwise} spine decomposition for supercritical super-Brownian motion with general branching mechanism 
(cf. \cite{KLMSR}) by studying $L^p$ convergence of  naturally underlying additive martingales in the spirit of analogous arguments for branching 
particle diffusions due to Harris and Hardy \cite{HH2009}.  Amongst other ingredients, the Dynkin-Kuznetsov $\mathbb{N}$-measure plays a pivotal role 
in the analysis.

\bigskip


\bigskip

\end{abstract}

\vspace{0.5cm}

\section{Introduction}

Suppose that $X=\{X_t: t\geq 0\}$ is a (one-dimensional) $\psi$-super-Brownian motion with general branching mechanism $\psi$ taking the form
\begin{equation}
\psi(\lambda) = -\alpha \lambda + \beta\lambda^2 + \int_{(0,\infty )} ({\rm e}^{-\lambda x} - 1 + \lambda x )\nu({\rm d}x),
\label{branch-mech}
\end{equation}
for $\lambda \geq 0$ where $\alpha = - \psi'(0+)\in(0,\infty)$, $\beta\geq 0$ and
$\nu$ is a measure concentrated on $(0,\infty)$ which satisfies
$\int_{(0,\infty)}(x\wedge x^2)\nu({\rm d}x)<\infty$. 
Let
$\mathcal{M}_F(\mathbb{R})$ be  the space of finite measures on $\mathbb{R}$
and note that $X$ is a $\mathcal{M}_F(\mathbb{R})$-valued Markov
process under $\mathbb{P}_\mu$ for each $\mu\in\mathcal{M}_F(\mathbb R)$, where $\mathbb{P}_\mu$ is \textcolor{black}{the} law of $X$ with initial 
configuration $\mu$. 
We shall
 use standard inner product notation, for
$f\in C_b^+(\mathbb{R})$, the space of positive, uniformly bounded, continuous functions on $\mathbb{R}$,  and $\mu\in\mathcal{M}_F(\mathbb{R})$,
\[
\langle f , \mu\rangle =  \int_{\mathbb{R} }f(x)\mu({\rm d}x).
\]
Accordingly we shall write $||\mu|| = \langle 1,\mu \rangle$.
Recall that the total mass of the process $X$, $\{||X_t||: t\geq 0\}$ is a continuous-state branching process with branching mechanism $\psi$. 
Such processes may exhibit explosive behaviour, however, under the conditions assumed above, $||X||$ remains finite at all  times.
We insist moreover that $\psi(\infty)=\infty$ which means that with positive probability the event $\lim_{t\uparrow\infty}||X_t||=0$ will occur. Equivalently this means that the total mass process does not have monotone increasing paths; see for example the summary in Chapter 10 of Kyprianou \cite{K}. 
The existence of these superprocesses processes is guaranteed by \cite{D, Dyn1993, Dyn2002}.

The following  standard result from the theory of superprocesses  describes the evolution of $X$ as a Markov process.
For all $f\in C_b^+(\mathbb{R})$  and $\mu\in\mathcal{M}_F(\mathbb{R})$,
\begin{equation}
 -\log\mathbb{E}_\mu({\rm e}^{- \langle f, X_t\rangle}) =  \int_{\mathbb{R} }u_f(x, t)\mu({\rm d}x), \, \mu\in\mathcal{M}_F(\mathbb{R}), \, t\geq 0,
 \label{prePDE}
\end{equation}
where  $u_f(x,t)$ is the unique positive solution to the evolution equation for $x\in\mathbb{R}$ and $t>0$
      \begin{equation}\label{PDE}
       \dfrac{\partial }{\partial
       t}u_f(x,t)=\dfrac{1}{2}\dfrac{\partial^2 }{\partial
       x^2}u_f(x,t)-\psi(u_f(x,t)),
       \end{equation}
with initial condition $u_f(x,0) = f(x)$. The reader is referred to Theorem 1.1 of  Dynkin \cite{Dyn1991}, Proposition 2.3 of Fitzsimmons \cite{Fitz} and Proposition 2.2 of Watanabe \cite{watanabe1968} for further details; see also Dynkin  \cite{Dyn1993, Dyn2002} and Engl\"ander and Pinsky \cite{EP} for a general overview.

Associated to the process $X$ is the following martingale  $Z(\lambda)=\{Z_t(\lambda),t\geq0\}$, where 
\begin{equation}\label{change1.1}
 Z_t(\lambda):={\rm e}^{\lambda c_\lambda t}\l {\rm e}^{\lambda\cdot} ,X_t\r,\,t\geq0,
\end{equation}
where $c_\lambda=\psi'(0+)/\lambda-\lambda/2$ and $\lambda\in\R$ (cf. \cite{KLMSR} Lemma 2.2). To see why this is a martingale note the following steps.
 Define for each $x\in\mathbb{R}$, $g\in C^+_b(\mathbb{R})$ and $\theta, t\geq 0$, $u_g^\theta(x,t) = - \log \mathbb{E}_{\delta_x}({\rm e}^{-\theta\langle g, X_t\rangle})$. With limits understood as $\theta\downarrow 0$, we  have $u_g(x,t)|_{\theta =0} =0$, moreover,    $v_g(x, t): =\mathbb{E}_{\delta_x}(\langle g, X_t\rangle) = \partial u_g^\theta(x,t)/\partial\theta |_{\theta = 0}$. Differentiating in $\theta$ in (\ref{PDE}) shows that $v_g$ solves the equation
   \begin{equation}\label{PDE2}
       \dfrac{\partial }{\partial
       t}v_g(x,t)=\dfrac{1}{2}\dfrac{\partial^2 }{\partial
       x^2}v_g(x,t)-\psi'(0+)v_g(x,t),
       \end{equation}
       with $v_g(x,0) = g(x)$.  Classical Feynman-Kac theory tells us that (\ref{PDE2}) has a unique solution which is necessarily equal to 
$\Pi_x({\rm e}^{-\psi'(0+) t} g(\xi_t))$ where $\{\xi_t : t\geq 0\}$ is a Brownian motion issued from $x\in\mathbb{R}$  under the measure $\Pi_x$.
The above procedure also works for $g(x) = \textcolor{black}{ {\rm e}^{\lambda x}}$ in which
case we easily conclude that for all $x\in\mathbb{R}$ and $t\geq 0$,
$\textcolor{black}{{\rm e}^{\lambda c_\lambda t}}\mathbb{E}_{\delta_x}(\langle \textcolor{black}{{\rm e}^{\lambda
\cdot}}, X_t \rangle) = \textcolor{black}{{\rm e}^{\lambda x}}$. Finally, the martingale
property follows using the previous equality together with the
Markov branching property associated with $X$. 
Note that as a positive martingale, it is automatic that 
$$\lim_{t\uparrow\infty}Z_t(\lambda)=Z_\infty(\lambda)$$
$\mathbb{P}_{\mu}$-almost surely for all $\mu\in\mathcal{M}_F(\mathbb{R})$ such that $\langle {\rm e}^{\lambda \cdot}, \mu\rangle<\infty$.

\bigskip

The purpose  of this note is to demonstrate the robustness of a new path decomposition of our $\psi$-super-Brownian motion by studying the $L^p$-convergence of the martingales $Z(\lambda)$. 
Specifically we shall prove the following theorem.

\begin{theorem} \label{Lptheorem} Assume that $p\in(1,2]$,  $\int_{(0,\infty)}r^p\nu({\rm d}r)< \infty$ and
$p\lambda^2<-2\psi^\prime(0+)$. Then
$Z_t(\lambda)$  converges to $Z_\infty(\lambda)$ in $L^p(\P_{\mu})$, for
all $\mu\in\mathcal{M}_F(\mathbb{R})$
such that $\langle {\rm e}^{\lambda\cdot}, \mu\rangle$ and $\l {\rm
e}^{\lambda p\cdot },\mu\r$ are finite.
\end{theorem}

The method of proof we use is quite similar to the one used in Harris and Hardy \cite{HH2009} for branching Brownian motion, where a pathwise spine decomposition functions as the key instrument of the proof. Roughly speaking, in that setting, the spine decomposition says that under a change of measure, the law of the branching Brownian motion has the same law as an immortal Brownian diffusion (with {\rm d}rift) along the path of which independent copies of the original branching Brownian motion immigrate at times which form a Poisson process. Until recently such a spine decomposition for superdiffusions was only available in the literature in a weak form; meaning that it takes the form of a semi-group decomposition. See the original paper of Evans \cite{Evans} as well as, for example amongst others, Engl\"ander and Kyprianou \cite{EK}.
Recently however Kyprianou et al. \cite{KLMSR} give a pathwise spine decomposition which provides a natural analogue to the pathwise spine decomposition for branching Brownian motion. Amongst other ingredients, the Dynkin-Kuznetsov $\mathbb{N}$-measure plays a pivotal \textcolor{black}{role} in describing the immigration off an immortal particle. We give a description of this new spine decomposition in the next section and thereafter we proceed to the proof of Theorem \ref{Lptheorem} in Section 3.

\section{Spine decomposition}
 For each $\lambda\in\R$ and $\mu\in\mathcal{M}(\mathbb{R})$ satisfying $\l e^{\lambda \cdot}, \mu\r $, we introduce the following martingale change of measure
\begin{equation}\label{change1}
 \frac{{\rm d} \P_{\mu}^{\lambda}}{{\rm d}\P_{\mu}}\bigg|_{\F_t}=\frac{Z_t(\lambda)}{\l e^{\lambda \cdot}, \mu\r},\,t\geq0,
\end{equation}
where $\F_t:=\sigma(X_s,s\leq t)$. 
The preceding change of measure induces the  {\it spine decomposition} of $X$ alluded to above. To describe it in detail we need some more ingredients.

According to Dynkin and Kuznetsov \cite{Dynkin-Kuznetsov} there exists a collection of measures $\{\N_x,x\in\R\}$, defined on the same probability space as $X$, such that 
\begin{equation}\label{N-measure}
\N_x\left(1-{\rm e}^{-\langle f,X_t\rangle}\right)=u_f(x,t),\,\,x\in\R, t\geq 0.
\end{equation}
Roughly speaking, the branching property tells us that for each $n\in\mathbb{N}$,  the measures $\mathbb{P}_{\delta_x}$ can be written as the $n$-fold convolution of $\mathbb{P}_{\frac{1}{n}\delta_{x}}$ which indicates that, on the trajectory space of the superprocess, $\mathbb{P}_x$ is infinitely divisible. Hence the role of $\mathbb{N}_x$ in (\ref{N-measure}) is analogous to that of the L\'evy measure for positive real-valued random variables.

From identity 
(\ref{N-measure}) and equation (\ref{prePDE}), it is straightforward to deduce that
\begin{equation}\label{N-measure-mean}
\mathbb{N}_x(\langle f, X_t \rangle) = \mathbb{E}_{\delta_x}(\langle f, X_t \rangle),
\end{equation}
whenever $f\in C_b^+(\mathbb{R})$.

For each $x\in\R$, let $\Pi_x$ be the law of a Brownian motion $\xi:=\{\xi_t : t\geq 0\}$ issued from $x$.
If $\Pi_x^{\lambda}$ is the law under which $\xi$ is a Brownian  motion with {\rm d}rift $\lambda\in\mathbb{R}$ and issued from $x\in\mathbb{R}$, 
 then for each $t\geq 0$,
\begin{equation}
  \frac{{\rm d}\Pi_x^{\lambda}}{{\rm d}\Pi_x}\bigg|_{\mathcal{G}_t}={\rm e}^{\lambda(\xi_t-x)-\frac{1}{2}\lambda^2t},\,t\geq0,
  \label{girsanov}
\end{equation}
where $\mathcal{G}_t:=\sigma(\xi_s,s\leq t)$. 
For convenience we shall also introduce the measure 
\begin{equation}
\Pi^\lambda_\mu(\cdot) : = \frac{1}{\langle {\rm e}^{\lambda \cdot}, \mu\rangle}\int {\rm e}^{\lambda x}\mu({\rm d}x)\Pi_x^\lambda(\cdot),
\label{randomization}
\end{equation}
for all $\lambda\in\mathbb{R}$. In other words, $\Pi^\lambda_\mu$ has the law of a Brownian motion with drift at rate $\lambda$ with an initial position which has been independently randomised in a way that is determined by $\mu$.

Now fix $\mu\in\mathcal{M}_F(\mathbb{R})$ and $x\in\mathbb{R}$ and let us define a measure-valued process $\Lambda:=\{\Lambda_t,t\geq0\}$ as follows:
\begin{enumerate}[(iii)]
\item[(i)] Take a copy of the process $\xi=\{\xi_t,t\geq0\}$ under $\Pi_x^{\lambda}$, we shall refer to this process as the {\it spine}.
\item[(ii)] Suppose that $\mathbf{n}$ is
                  a Poisson point process such that, for $t\geq 0$, given the
                  spine $\xi$, ${\bf n}$ issues superprocess $X^{{\bf n},t}$ at space-time position $(\xi_t, t)$ with rate $ {\rm d}t\times 2 \beta{\rm d}\mathbb{N}_{\xi_t}$.

\item[(iii)] Suppose that $\mathbf{m}$ is a
            Poisson point process  such that, independently of ${\bf n}$, given the spine $\xi$,  $\mathbf{m}$ issues a superprocess $X^{{\bf m}, t}$ at space-time point $(\xi_t, t)$ with initial mass $r$ at rate
            ${\rm d}t\times r \nu({\rm d}r)\times {\rm d} \mathbb{P}_{r\delta_{\xi_t}}$.                   
               
\end{enumerate} 

Note in particular that, when $\beta >0$, the rate of immigration under the process $\mathbf{n}$ is infinite and moreover, each process that immigrates is issued with zero mass. One may therefore think of $\mathbf{n}$ as a process of continuous immigration. In contrast, when $\nu$ is a non-zero measure, processes that immigrate under $\mathbf{m}$ have strictly positive initial mass and therefore contribute to path discontinuities of $||X||$.

Now, for each $t\geq0$, we define 
\begin{equation}\label{Lambda}
\Lambda_t = {X}'_t  +  X_t^{(\mathbf{n})} +  X_t^{(\mathbf{m})},
\end{equation}
where $\{X'_t: t\geq 0\}$ is an independent copy of $(X, \mathbb{P}_{\color{black}{\mu}})$,
\begin{equation*}
 X_t^{(\mathbf{n})}=\sum_{s\leq t:\mathbf n}X^{\mathbf n,
 s}_{t-s},\, t\geq 0 \qquad
\mbox{and}\qquad
 X_t^{(\mathbf{m})}=\sum_{s\leq t:\mathbf{m}}X^{\mathbf m, s}_{t-s},\, t\geq 0.
\end{equation*}
In the last two equalities we understand the first sum to be over times for which  ${\bf n}$ experiences points and the second sum is 
understood similarly. Note that since the processes $ X^{(\mathbf{n})}$ and $ X^{(\mathbf{m})}$ are initially zero valued it is clear that since $X'_0 = \mu$ then $\Lambda_0 = \mu$. In that case, we use the notation $\widetilde{\mathbb{P}}^{\lambda}_{\mu,x}$ to denote the law of the pair $( \Lambda, \xi)$. Note also that the pair $(\Lambda, \xi)$ \textcolor{black}{is} a time-homogenous Markov process. We are interested in the case that the initial position of the spine $\xi$ is randomised using the measure $\mu$ via (\ref{randomization}). In that case we shall write
\[
\widetilde{\mathbb{P}}^{\lambda}_\mu(\cdot) = \frac{1}{\l {\rm e}^{\lambda \cdot}, \mu\r }\int_{\mathbb{R}} {\rm e}^{\lambda x}\mu({\rm d}x)\widetilde{\mathbb{P}}^{\lambda}_{\mu, x}(\cdot)
\]
 for short. The next theorem identifies the process $\Lambda$ as the {\it pathwise} spine decomposition of $(X,\mathbb{P}^{\lambda}_{\mu})$ and in particular it shows that as a process on its own $\Lambda$ is Markovian.
\begin{teo}(Theorem 5.1, \cite{KLMSR})\label{spine-decomp} For all $\mu\in\mathcal{M}_F(\mathbb{R})$ such that $\langle {\rm e}^{\lambda\cdot}, \mu\rangle<\infty$, $(X, \mathbb{P}^{\lambda}_{\mu})$ and 
$(\Lambda, \widetilde{\mathbb{P}}^{\lambda}_{\mu})$ are equal in law.
\end{teo}

\section{Proof of Theorem \ref{Lptheorem}}
From the last section we have the following spine decomposition of the martingale (\ref{change1.1}),
\begin{equation}\label{spine}
 Z_t^\Lambda(\lambda)=Z_t^\prime(\lambda)+\sum_{s\leq t:{\bf n}}{\rm e}^{\lambda c_\lambda s}Z^{{\bf n},s}_{t-s}(\lambda)
+\sum_{s\leq t:{\bf m}}{\rm e}^{\lambda c_\lambda s}Z^{{\bf m},s}_{t-s}(\lambda),
\end{equation}
where $Z^\prime(\lambda)$ is an independent copy of $Z(\lambda)$ under $\P_{\mu}$, 
\begin{equation*}
Z_{t-s}^{\mathbf{n},s}:={\rm e}^{\lambda c_\lambda(t-s)}\langle {\rm e}^{\lambda\cdot},X_{t-s}^{\mathbf{n},s}\rangle,
\end{equation*}
and 
\begin{equation*}
Z_{t-s}^{\mathbf{m},s}:={\rm e}^{\lambda c_\lambda(t-s)}\langle {\rm e}^{\lambda\cdot},X_{t-s}^{\mathbf{m},s}\rangle.
\end{equation*}


Since $\{Z_t(\lambda),t\geq0\}$ is a martingale and we assume that  $p\in(1,2]$, then Doob's submartingale inequality
 tells us that  $Z(\lambda)$ converges in $L^p(\P_{\mu})$ as soon as we can show that $\sup_{t\geq 0}\E_{\mu}(Z_t(\lambda)^p)<\infty$. To this end, and with the above pathwise spine decomposition in hand we may now proceed to ad{\rm d}ress the analogue of the proof for branching Brownian motion given in \cite{HH2009}.

First note that, for all $p\in(1,2]$,
\begin{equation}
 \E_{\mu}(Z_t(\lambda)^p)=\E_{\mu}^\lambda(Z_t(\lambda)^q) =\widetilde{\E}_{\mu}^\lambda(Z^\Lambda_t(\lambda)^q),\,\,\mbox{for all}\,\,t\geq0,
 \label{supfin}
\end{equation}
where $q=p-1$. By Jensen's inequality we have that, for all $q\in(0,1]$
\begin{eqnarray}\label{cond-Jensen}
\nonumber\lefteqn{\widetilde{\E}_{\mu}^\lambda\left(Z^\Lambda_t(\lambda)^q\mid\xi\right)}\\&\leq&\nonumber\left[\widetilde{\E}_{\mu}^\lambda\left(Z^\Lambda_t(\lambda)\mid\xi\right)\right]^q
\\&\leq&\l{\rm e}^{\lambda \cdot}, \mu\r^{\color{black}{q}}+\left[\widetilde{\E}_{\mu}^\lambda\left(\sum_{s\leq t:{\bf n}}{\rm e}^{\lambda c_\lambda s}Z^{{\bf n},s}_{t-s}(\lambda)\bigg|\xi\right)\right]^q+\left[\widetilde{\E}_{\mu}^\lambda\left(\sum_{s\leq t:{\bf m}}{\rm e}^{\lambda c_\lambda s}Z^{{\bf m},s}_{t-s}(\lambda)\bigg|\xi\right)\right]^q,
\end{eqnarray}
\textcolor{black}{to get the last inequality we have used the fact that $(\sum_iu_i)^q\leq\sum_iu_i^q$ with $u_i\geq0$}. On the one hand, recalling from  (\ref{N-measure-mean}) that $\N_{\xi_s}[Z_{t-s}(\lambda)]=\E_{\xi_s}[Z_{t-s}(\lambda)]$, we obtain
\begin{eqnarray}\label{cond-Jensen.uno}
\widetilde{\E}_{\mu}^\lambda\left(\sum_{s\leq t:{\bf n}}{\rm e}^{\lambda c_\lambda s}Z^{{\bf n},s}_{t-s}(\lambda)\bigg|\xi\right)&=&\nonumber\int_0^t{\rm e}^{\lambda c_\lambda s}\N_{\xi_s}[Z_{t-s}(\lambda)]{\rm d}s
\\&=&\int_0^t{\rm e}^{\lambda (\xi_s+c_\lambda s)}{\rm d}s.
\end{eqnarray}
On the other hand, we have that 
\begin{eqnarray}\label{cond-Jensen.dos}
\nonumber\widetilde{\E}_{\mu}^\lambda\left(\sum_{s\leq t:{\bf m}}{\rm e}^{\lambda c_\lambda s}Z^{{\bf m},s}_{t-s}(\lambda)\bigg|\xi\right)&=&\widetilde{\E}_{\mu}^\lambda\left[\widetilde{\E}_{\mu}^\lambda\left(\sum_{s\leq t:{\bf m}}{\rm e}^{\lambda c_\lambda s}Z^{{\bf m},s}_{t-s}(\lambda)\bigg|\xi,\mathbf{m}\right)\bigg|\xi\right]
\\&=&\nonumber\widetilde{\E}_{\mu}^\lambda\left(\sum_{s\leq t:{\bf m}}m_s{\rm e}^{\lambda (\xi_s+c_\lambda s)}\bigg|\xi\right)
\\&=&\sum_{s\leq t:{\bf m}}m_s{\rm e}^{\lambda (\xi_s+c_\lambda s)},
\end{eqnarray}
where for $s\geq 0$, $m_s = ||X^{\mathbf{m}, s}_0||$. In particular note that the process $\{m_t : t\geq 0\}$ is a Poisson point process on $(0,\infty)^2$, independent of $\xi$, with intensity ${\rm d} t \times r\nu({\rm d} r)$.
Then,  putting (\ref{cond-Jensen.uno}) and (\ref{cond-Jensen.dos}) into (\ref{cond-Jensen}), making use again of the inequality $(\sum_i u_i)^q \leq \sum_i u_i^q$ where $u_i\geq 0$ for all $i$, we obtain
\begin{eqnarray}\label{cond-Jensen1}
\nonumber\widetilde{\E}_{\mu}^\lambda\left(Z^\Lambda_t(\lambda)^q\mid\xi\right)&\leq&\l{\rm e}^{\lambda \cdot}, \mu\r^{\color{black}{q}}+\left(\int_0^t{\rm e}^{\lambda (\xi_s+c_\lambda s)}{\rm d}s\right)^q+\sum_{s\leq t:{\bf m}}m_s^q{\rm e}^{q\lambda(\xi_s+c_\lambda s)}
\\&\leq&\l{\rm e}^{\lambda \cdot}, \mu\r^{\color{black}{q}}+\left(\int_0^\infty {\rm e}^{\lambda (\xi_s+c_\lambda s)}{\rm d}s\right)^q+\sum_{s\geq0:{\bf m}}m_s^q{\rm e}^{q\lambda(\xi_s+c_\lambda s)}.
\end{eqnarray}
Taking expectations again in  (\ref{cond-Jensen1}) gives us that, for all $t\geq0$,
\begin{equation}\label{q-moment}
\widetilde{\E}_{\mu}^\lambda(Z^\Lambda_t(\lambda)^q)\leq {\color{black} \langle {\rm e}^{\lambda\cdot},\mu\rangle^q}+\Pi_\mu^\lambda\left[\left(\int_0^\infty {\rm e}^{\lambda (\xi_s+c_\lambda s)}{\rm d}s\right)^q\right]+\widetilde{\E}_{\mu}^\lambda\left(\sum_{s\geq0:{\bf m}}m_s^q{\rm e}^{q\lambda(\xi_s+c_\lambda s)}\right).
\end{equation}
We know that, under $\Pi^\lambda_\mu$, the process $\xi$ is a Brownian motion with {\rm d}rift $\lambda$. Thus, with respect to the same measure,  $\xi_s+c_\lambda s$ is a Brownian motion with
drift $\lambda+c_\lambda$ which is strictly negative for $\lambda\in(0,\sqrt{-2\psi^\prime(0+)})$. Note that this latter condition in particular holds under assumption that $p\lambda^2<-2\psi^\prime(0+)$ and $p>1$. From Section 2 of Maulik and Zwart \cite{MZ} we can conclude that 
 \begin{equation*}
 \Pi^\lambda_0\left(\int_0^\infty {\rm e}^{\lambda (\xi_s+c_\lambda s)}{\rm d}s\right)<\infty,
\end{equation*}
which in turn implies that, for all $q\in(0,1]$,
 \begin{equation*}
 \Pi^\lambda_0\left[\left(\int_0^\infty {\rm e}^{\lambda (\xi_s+c_\lambda s)}{\rm d}s\right)^q\right]<\infty,
\end{equation*}
and hence
 \begin{eqnarray}
 \Pi^\lambda_\mu\left[\left(\int_0^\infty {\rm e}^{\lambda (\xi_s+c_\lambda s)}{\rm d}s\right)^q\right]& =& \frac{1}{\l {\rm e}^{\lambda\cdot},\mu\r} \int {\rm e}^{\lambda x }\mu({\rm d} x) \Pi_0^\lambda\left[\left(\int_0^\infty {\rm e}^{\lambda (x+\xi_s+c_\lambda s)}{\rm d}s\right)^q\right] \notag \\
 &=& \frac{\l {\rm e}^{\lambda p \cdot}, \mu\r}{\l {\rm e}^{\lambda \cdot},\mu \r } \Pi^\lambda_0\left[\left(\int_0^\infty {\rm e}^{\lambda (\xi_s+c_\lambda s)}{\rm d}s\right)^q\right]<\infty.
\label{N-term}
\end{eqnarray}
It remains to prove that the last term in (\ref{q-moment}) is finite. This can be done by computing  the expectation directly. We obtain, 
\begin{eqnarray*}
\widetilde{\E}_{\mu}^\lambda\left(\sum_{s\geq0:{\bf m}}m_s^q{\rm e}^{q\lambda(\xi_s+c_\lambda s)}\right)&=&\int_0^\infty {\rm d}s\int_0^\infty r\nu({\rm d}r)r^q\Pi_\mu^\lambda\left({\rm e}^{q\lambda(\xi_s+c_\lambda s)}\right)
\\
&=&\int_0^\infty {\rm d}s\int_0^\infty r^p\nu({\rm d}r) \frac{1}{\l {\rm e}^{\lambda\cdot},\mu\r} \int {\rm e}^{\lambda x }\mu({\rm d} x) \Pi_0^\lambda\left({\rm e}^{q\lambda(x+\xi_s+c_\lambda s)}\right)\\
&=&{\rm e}^{q\lambda x} \frac{\l {\rm e}^{\lambda p \cdot}, \mu\r}{\l {\rm e}^{\lambda \cdot},\mu \r }\int_0^\infty r^p\nu({\rm d}r)\int_0^\infty \Pi_0^\lambda\left({\rm e}^{q\lambda(\xi_s+c_\lambda s)}\right){\rm d}s.
\end{eqnarray*}
Note that, 
\begin{eqnarray*}
\Pi_0^\lambda\left({\rm e}^{q\lambda(\xi_s+c_\lambda s)}\right)&=&\exp\{ qs\lambda^2+s(q\lambda)^2/2 + qs\psi'(0+)-qs\lambda^2/2\}\\
&=&\exp\{ qs( p\lambda^2/2+ \psi'(0+))\}
\end{eqnarray*}
for all $s\geq0$. Moreover, this expectation has a negative exponent as soon as $p\lambda^2<-\psi^\prime(0+)$. Together 
with the assumption $\int_0^\infty r^p\nu({\rm d}r)<\infty$ we conclude that
\begin{equation}\label{P-term}
\widetilde{\E}_{\mu}^\lambda\left(\sum_{s\geq0:{\bf m}}m_s^q{\rm e}^{q\lambda(\xi_s+c_\lambda s)}\right)<\infty.
\end{equation} 
Finally, from (\ref{q-moment})-(\ref{P-term}) we get that
\begin{equation*}
\sup_{t\geq 0}\widetilde{\E}_\mu^\lambda\left(Z^\Lambda_t(\lambda)^q\right)<\infty,
\end{equation*}
which, in combination with (\ref{supfin}), completes the proof.
 \ter


\section*{Acknowledgments} 
The second author would like to thank the University of Bath, where most of this research was done. He also acknowledges the financial support of
 CONACYT-Mexico grant number 129076.

\noindent{\sc Department of Mathematical Sciences, University of Bath, Claverton Down, Bath BA2 7AY, United Kingdom.} \\{\it E-mail:} {a.kyprianou@bath.ac.uk}

\noindent{\sc Departamento de Matem\'aticas, Universidad de Guanajuato,
Jalisco S/N Mineral de Valenciana, Guanajuato, Gto.\ C.P.\ 36240,
M\'exico.} \\{\it E-mail:} {amurillos@ugto.mx}


\begin{thebibliography}{99}




\bibitem{D} E.B. Dynkin (1991):
         Branching particle systems and superprocesses, {\em  Ann. Probab.} {\bf 19(3)}, 1157--1194.

\bibitem{Dyn1991} E.B. Dynkin (1991): A probabilistic approach to one class of non-linear differential equations. {\it Probab. Th. Rel. Fields}. {\bf 89}, 89--115.

\bibitem{Dyn1993} E.B. Dynkin (1993): Superprocesses and Partial Differential Equations. {\it Ann. Probab.} {\bf  21}, 1185--1262.


\bibitem{Dyn2002} E.B. Dynkin (2002): {\it Diffusions, Superdiffusions and Partial Differential Equations.} AMS, Providencem R.I.



\bibitem{Dynkin-Kuznetsov} E.B. Dynkin and S.E. Kuznetsov (2004): $\mathbb{N}$-measures for branching Markov exit systems and their applications to 
differential equations.  {\it Probab. Th. Rel. Fields}. {\bf 130}, 135--150.


\bibitem{EP} Engl\"ander, J. and Pinsky, R. (1999): On the construction and support properties of measure-valued diffusions on $D\subseteq R\sp d$ with spatially dependent branching.
{\it Ann. Probab.} {\bf 27}, 684--730.




       \bibitem{Evans} S.N. Evans (1993):  Two representations of a superprocess. {\it Proc. Royal. Soc. Edin.} {\bf 123A}, 959--971.

   \bibitem{Fitz}  P.J. Fitzsimmons (1988): Construction and regularity of measure-valued Markov branching processes. {\it Israeli J. Math.} {\bf 64}, 337--361.


\bibitem{EK} J.  Engl\"{a}nder  and A.E. Kyprianou (2004): Local
extinction
versus local exponential growth for spatial
branching processes. {\it  Ann. Probab.} {\bf  32}, 78--99.

\bibitem{HH2009} R. Hardy and S.C. Harris (2009): A spine approach to branching diffusions with applications to $\mathcal{L}^p$-convergence of martingales. 
{\it S\'eminaire de Probabilit\'es}, {\bf XLII}, 281-330. 




              
              \bibitem{Kyp} Kyprianou, A.E. (2004): Travelling wave solutions to the K-P-P equation: alternitives to Simon Harris' probabilistic analysis.
{\it Ann. Inst. H. Poincar\'e}.  {\bf 40}, 53--72. 
               
              \bibitem{K} Kyprianou, A.E. (2006): {\it Introductory lectures on fluctuations of L\'evy processes with applications.} Springer.

\bibitem{KLMSR} A.E. Kyprianou, R.-L. Liu, A. Murillo-Salas and Y.-X. Ren. (2011): Supercritical super-Brownian motion with a general branching mechanism and travelling waves. {\it To appear in Ann. Inst. H. Poincar\'e}.

\bibitem{MZ} K. Maulik and B. Zwart (2006): Tail asymptotics for exponential functionals of L\'evy processes. {\it Stoch. Proc. Appl.} {\bf 116}, 156--177.




\bibitem{watanabe1968} S. Watanabe (1968): A limit theorem of branching processes and continuous-state branching processes. {\em J. Math. Kyoto Univ.} {\bf 8}, 141--167.







\end{thebibliography}
\end{document}